\documentclass[12pt]{article}
\usepackage{german}
\usepackage{amsmath}
\usepackage{amssymb}
\usepackage{latexsym}

\newcommand{\ang}{\raisebox{0.2ex}{\scriptsize$\triangleright$}}

\newcommand{\mn}{\medskip\noindent}
\newcommand{\bn}{\bigskip\noindent}
\newcommand{\sn}{\smallskip\noindent}

\newcommand{\D}{{\cal{D}}}  
\newcommand{\cL}{{\cal{L}}}
\newcommand{\E}{{\cal{E}}}  
\newcommand{\K}{{\cal{K}}}  
\newcommand{\C}{{\cal{C}}}  
\newcommand{\B}{{\cal{B}}}  
\newcommand{\G}{{\cal{G}}}  
\newcommand{\A}{{\cal{A}}}  
\newcommand{\N}{{\cal{N}}}  
\newcommand{\cP}{{\cal{P}}} 
\newcommand{\Q}{{\cal{Q}}}  
\newcommand{\X}{{\cal{X}}}  
\newcommand{\sS}{{\cal{S}}} 
\newcommand{\Hh}{{\cal{H}}}  
\newcommand{\cO}{{\cal{O}}}

\newcommand{\g}{{\mathfrak{g}}}

\newcommand{\pb}{{\mathsf{p}}}
\newcommand{\bx}{{\mathsf{x}}}
\newcommand{\by}{{\mathsf{y}}}

\newcommand{\bbbr}{{\mathbb{R}}}
\newcommand{\bbbc}{{\mathbb{C}}}
\newcommand{\bbbb}{{\mathbb{B}}}

\begin{document}
\begin{center} 
{\LARGE On Well-behaved Unbounded Representations of $\ast$-Algebras}
\end{center}

\centerline{Konrad Schm\"udgen}
\centerline{Fakult\"at f\"ur Mathematik und Informatik}
\centerline{Universit\"at Leipzig, Augustusplatz 10, 04109 Leipzig, Germany}
\centerline{E-mail: schmuedg@mathematik.uni-leipzig.de}

\renewcommand{\baselinestretch}{1.0}

\mn
{\bf Abstract:} A general approach to the well-behaved unbounded 
$\ast$-representations of a $\ast$-algebra $\X$ is proposed. Let $\B$ be 
a normed $\ast$-algebra equipped with a left action $\ang$ of 
$\X$ on $\B$ such that
 $(x\ang a)^+ b=a^+(x^+\ang b)$ for $a,b\in\B$ and $x\in\X$. Then the pair 
 $(\X,\B)$ is called a {\it compatible pair}. For any 
continuous non-degenerate $\ast$-representation $\rho$ of 
$\B$ there exists a closed $\ast$-representation 
$\rho^\prime$ of $\X$ such that $\rho^\prime(x)\rho(b)=\rho(x\ang b)$, 
where $x\in\X$ and $b\in\B$. The $\ast$-representations $\rho^\prime$ 
are called the well-behaved $\ast$-representations associated with
the compatible pair $(\X,\B)$. A number of examples are developed in detail.


\mn
{\bf 0. Introduction}

\sn
Unbounded representations of general $\ast$-algebras in Hilbert space occur 
in various branches of mathematics and mathematical physics such as 
representation theory of Lie algebras, algebraic quantum field theory, 
the theory of quantum groups and quantum algebras. One of the natural 
questions is to ask for a description or classification of {\it all} 
$\ast$-representations of the corresponding $\ast$-algebra. But it turns 
out that this is not a well-posed problem for general $\ast$-algebras. 
In order to explain this, let $\X$ be the 
$\ast$-algebra $\bbbc[x,y]$ of all polynomials in two commuting hermitean 
indeterminants $x$ and $y$ or the $\ast$-algebra $A(p,x)$ of two 
hermitean generators $p$ and $x$ satisfying the Heisenberg commutation 
relation $px-xp=-i$. In both cases it seems to be impossible to classify 
in a reasonable way all $\ast$-representations of $\X$  even if we assume that 
the images of the generators $x,y$ resp. $p,x$ are essentially self-adjoint 
(see, for instance, [S1] Chapter 9, and [S2]). The large 
variety of such $\ast$-representations of the 
$\ast$-algebra $\bbbc[x,y]$ is illustrated by the following 
result ([S1], Theorem 9.4.1): For any properly infinite von Neumann 
algebra $\N$ on a separable Hilbert space there exists a 
$\ast$-representation $\rho$ of the polynomial algebra $\bbbc[x,y]$ such 
that the operators $\rho(x)$ and $\rho(y)$ are essentially self-adjoint 
and the spectral projections of these operators generate the von Neumann 
algebra $\N$.

In most situations it suffices to know some class of 
 ''nice'' $\ast$-representations 
of the $\ast$-algebra which is characterized by means of additional 
requirements in order to exclude pathological behaviour of operators. 
In what follows we shall call these 
$\ast$-representations {\it well-behaved}. In earlier papers 
(see [S1], Chapters 9 and 10, and [S3]) we have called them 
integrable representations because of the commonly used 
terminology in representation theory of Lie algebras. For the 
$\ast$-algebras $\bbbc[x,y]$ and $A(p,x)$ it is easy to guess how 
to define well-behaved $\ast$-representations. A $\ast$-representation 
$\rho$ of $\bbbc[x,y]$ is called well-behaved if $\rho$ is self-adjoint 
(see [P1], [P2] or [S1] for this notion) and if $\rho (x)$ and $\rho(y)$ 
are essentially self-adjoint 
operators such that their spectral projections commute. In the case 
$\X=A(p,q)$ the latter condition should be replaced by the requirement 
that $P:=\overline{\rho(p)}$ and $X:=\overline{\rho(y)}$ are self-adjoint 
operators satisfying the Weyl relation 
$e^{itP} e^{isX} = e^{ist} e^{isX} e^{itP}, s,t\in\bbbr$.

Many papers of the mathematical physics literature dealing with 
unbounded $\ast$-representations claim to determine 
{\it all} $\ast$-representations of the $\ast$-algebra. 
However, a closer look at the 
proofs shows that often hidden additional assumptions are used and that
only a particular class of well-behaved $\ast$-representations is  
investigated. For instance, for $\ast$-algebras related to the canonical 
commutation relations or for $q$-oscillator algebras it often 
required that a vacuum vector exists or that certain operators have a 
complete set of eigenvectors.

There is no general method to select the well-behaved $\ast$-representations 
of a given $\ast$-algebra. 
Also it is important to stress that the choice of well-behaved 
$\ast$-representations may depend on the aim of considerations. 
The examples developed in Section 4 show that for the same $\ast$-algebra  
there are various natural candidates for the definition of well-behavedness. 
The additional conditions selecting well-behaved $\ast$-representations 
depend, generally speaking, on the underlying $\ast$-algebra.  

In this paper we propose a general approach to the study of 
well-behaved $\ast$-representations. The idea is easily explained as 
follows: Let $\X$ be a $\ast$-algebra and let $\B$ be 
a normed $\ast$-algebra equipped with a left action, written $x\ang b$, of 
$\X$ on $\B$ satisfying the compatibility condition
 $(x\ang a)^+ b=a^+(x^+\ang b)$ for $a,b\in\B$ and $x\in\X$. We shall call 
 such a pair $(\X,\B)$ a {\it compatible pair}. Then, for any 
continuous non-degenerate $\ast$-representation $\rho$ of the 
normed $\ast$-algebra $\B$ there exists a closed $\ast$-representation 
$\rho^\prime$ of $\X$ such that $\rho^\prime(x)\rho(b)=\rho(x\ang b)$, 
where $x\in\X$ and $b\in\B$. The $\ast$-representations $\rho^\prime$ of 
$\X$ obtained in this way are called the well-behaved $\ast$-representations 
of $\X$ associated with the compatible pair $(\X,\B)$.

The main purpose of this paper is to develop a number of important 
examples and show that they fit into this context. In Section 2 we discuss 
some compatible pairs $(\X,\B)$ for the polynomial algebra 
$\X=\bbbc[x_1,{\dots},x_n]$. In Section 3 we treat that the $G$-integrable 
representations of the enveloping algebra $\E(\g)$ of the Lie algebra$\g$ 
of a Lie group $G$. Here $\B$ is the $\ast$-algebra 
$C^\infty_0(G)$ with convolution multiplication. In Section 4 we study 
various $\ast$-compatible pairs $(\X,\B)$ by using the Weyl calculus of 
pseudodifferential operators. Among others, various classes of 
well-behaved $\ast$-representations of the coordinate 
$\ast$-algebra $\X=\cO(\bbbr^2_q)$ of the real quantum plane are 
considered in this approach. In Section 5 we consider the quantum 
group $SU_q(1,1)$. The paper closes with a short outlook in Section 6.

Let us collect some definitions and facts on unbounded operator algebras 
and unbounded $\ast$-representations used in what follws. 
More details can be found in the monograph [S1]; see also [In] or [P1]. 

Let $\D$ be a dense linear subspace of a Hilbert space $\Hh$. Then the 
vector space
$$
\cL^+(\D)=\{x\in {\rm End}~\D :\D\subseteq\D(x^\ast), x^\ast\D\subseteq \D\}
$$ 
is a unital $\ast$-algebra with operator product as multiplication and 
the restriction $x^+:= x^\ast i\D$ to $\D$ of the adjoint operator 
$x^\ast$ to $\D$ as involution. The set $\bbbb(\D)$ of all bounded operator $b$ on $\Hh$ such that 
$b\Hh\subseteq\D$ and $b^\ast\Hh\subseteq\D$ is obviously a 
$\ast$-algebra. An $O^\ast$-algebra on the domainn $\D$ 
is a $\ast$-subalgebra of $\cL^+(\D)$ which contains the identity map 
of $\D$. By a $\ast$-representation of an abstract 
$\ast$-algebra $A$ (without unit in general) on a domain $\D$ 
we mean a $\ast$-homomorphism $\rho$ of $\A$ into the $\ast$-algebra 
$\cL^+(\D)$. The $\ast$-representation $\rho$ of $A$ is said to be closed if 
$\D$ is the intersection of all domains $\D(\overline{\rho(a)})$, $a\in A$, 
where the bar refers to the closure of the operator $\rho(a)$. The representation 
$\rho$ is called non-degenerate if
$\rho(A)\D$ is dense in the underlying Hilbert space $\Hh$. 
If $A$ is a normed $\ast$-algebra (that is, $\A$ is equipped with a
submultiplicative $\ast$-invariant norm $\parallel \cdot \parallel$), 
then a $\ast$-representation 
$\rho$ is called continuous if all operators are $\rho(a)$, $a \in A$, are
bounded and there exist a positive constant $C$ such that 
$\parallel \rho(a) \parallel  \leq C \parallel a \parallel$ for all $a\in A$, 
where $\parallel \rho(a) \parallel$ denotes the the operator norm of 
$\rho(a)$.  

\mn
{\bf 1. Compatible Pairs}

\sn
Let $\X$ be a $\ast$-algebra with unit element 1 and let $\B$ be a 
$\ast$-algebra (without unit in general). The involutions of $\X$ and 
$\B$ are denoted by $x\rightarrow x^+$ and $b\rightarrow b^+$, 
respectively. Suppose that the vector space $\B$ is a left $\X$-module 
with left action denoted by $\ang$, that is, there exists a linear 
mapping $\phi:\X\otimes \B\rightarrow \B$, written as 
$\phi (x\otimes b)= x\ang b$, such that $(xy)\ang b=x\ang (y\ang b)$ 
and $1\ang b=b$ for $x,y\in\X$ and $b\in\B$.

\mn
{\bf Proposition 1.} {\it Suppose that the left action of the $\ast$-algebra 
$\X$ on the $\ast$-algebra satisfies the condition 
\begin{equation}\label{cond}
(x\ang a)^+ b=a^+(x^+\ang b)~{\rm for~all}~x\in\X~{\rm and}~ a,b\in\B.
\end{equation}
Then, for any non-degenerate $\ast$-representation $\rho$ of the 
$\ast$-algebra $\B$ there exists a unique $\ast$-representation 
$\tilde{\rho}$ of the $\ast$-algebra $\X$ on the domain 
$\D(\tilde{\rho})=\rho (B)\D(\rho)$ such that
\begin{equation}\label{def}
\tilde{\rho}(x)(\rho(b)\varphi)=\rho(x\ang b)\varphi, x\in\X, b\in\B, 
\varphi\in\D(\rho).
\end{equation}
Let $\rho^\prime$ denote the closure of the $\ast$-representation 
$\tilde{\rho}$.}
{\bf Proof.} Let $\xi=\sum_i \phi (a_i)\varphi_i$ and $\eta=\sum_j\rho(b_j)\psi_j$ be vectors of the domain $\D(\tilde{\rho})$, where $a_ib_j\in\B$ and $\varphi_i, \psi_j\in\D(\rho)$. Let $x\in\X$. Using condition (\ref{cond}) and the assumption that $\rho$ is a $\ast$-representation of the $\ast$-algebra $\B$ we compute

\begin{align}\label{cor}
&~\quad\Big\langle {\sum\nolimits_i} \rho(x\ang a_i)\varphi_i, \eta\Big
   \rangle
  =\Big\langle {\sum\nolimits_i} \rho(x\ang a_i)\varphi_i,
   {\sum\nolimits_j}\pi(b_j)\psi_j\Big\rangle
  ={\sum\nolimits_{i,j}} \langle\varphi_i,\rho((x\ang a_i)^+b_j)\psi_j
   \rangle\nonumber\\
& ={\sum\nolimits_{i,j}} \langle \varphi_i,\rho(a_i^+(x^+\ang b_j))\psi_j
   \rangle 
  =\Big\langle {\sum\nolimits_i} \rho(a_i)\varphi_i,    {\sum\nolimits_j}\rho(x^+\ang b_j)\psi_j\Big\rangle
  =\Big\langle\zeta,{\sum\nolimits_j}\rho(x^+\ang b_j)\psi_j\Big\rangle
\end{align}
First we shall use relation (\ref{cor}) in order to show that equation (\ref{def}) defines unambignoulsy a linear operator $\tilde{\rho}(x)$ on the domain $\D(\tilde{\rho})$. In order to do so, it suffices to check that $\zeta\equiv\sum_i\rho(a_i)\varphi_i=0$ implies that $\sum_i\rho (x\ang a_i)\varphi_i=0$. Indeed, if $\zeta=0$, then 
it follows from (2) that $\langle \sum_i\rho(x\ang a_i)\varphi_i,\zeta\rangle=0$ for all $\eta\in\D(\tilde{\rho}).$ Since $\rho$ is non-degenerate, $\D(\tilde{\rho})$ is dense in the underlying Hilbert space and so $\sum_i\rho(x\ang a_i)\varphi_i=0$. Hence the operator $\tilde{\rho} (x)$ is well-defined.

{}From the properties of a left action if follows at once that $\tilde{\rho}$ 
is an algebra homomorphism of $\X$ into the linear operators acting on the 
domain $\D(\tilde{\rho})$ and leaving $\D(\tilde{\rho})$ invariant. 
In order to prove that $\tilde{\rho}$ preserves the involution, we 
combine equations (\ref{def}) and (\ref{cor}) and conclude that 
$\langle \tilde{\rho}(x)\zeta, \eta\rangle=\langle\zeta,
\tilde{\rho}(x^+)\eta\rangle$ for all $\zeta,\eta\in\D(\tilde{\rho})$. 
Thus, $\tilde{\rho}$ is indeed a $\ast$-representation of the $\ast$-algebra 
$\X$ on the domain $\D(\tilde{\rho})$.\hskip4.5cm$\Box$

\mn
{\bf Definition 2.} A {\it compatible pair} is a pair $(\X,\B)$ of a 
unital $\ast$-algebra $\X$ and a normed $\ast$-algebra $\B$ equipped with a 
left action of $\X$ on $\B$ satisfying condition  (\ref{cond}).
\mn

Our guiding example of a compatible pair is the following 

{\bf Example 3}: Let $\X$ be an $O^\ast$-algebra on a dense domain $\D$ 
of  a Hilbert space and let $\B$ be a $\ast$-subalgebra of $\bbbb(\D)$ 
such that $x b\in\B$ for all $x\in\X$ and $x\in\B$. We equipp the $\ast$-algebra $\B$ with the operator 
norm. Then there is a left action $\ang$ of $\X$ on $\B$ defined by the 
operator product, that is, $x\ang b:=bx,~x\in\X, b\in\B.$ 
It is not difficult to show 
that $(xb)^+ a= b^+ x^+ a$ for $a,b\in\B$ and $x\in\X$. Hence $(\X,\B)$ 
is a compatible pair. We call such a pair $(\X,\B)$ a 
{\it compatible $O^\ast$-pair} on the domain $\D$. In particular, 
$(\cL^+(\D),\bbbb(\D))$ is a compatible $O^\ast$-pair, because 
$xb\in\bbbb(\D)$ for $x\in\cL^+(\D)$ and $b\in\bbbb(\D)$. 

\mn
Now let $(\X,\B)$ be a compatible pair and let $\rho$ be a 
continuous $\ast$-representation of the normed $\ast$-algebra $\B$ on a 
Hilbert space $\Hh$. Then it is clear that $\tilde{\rho}(\X)$ is an 
$O^\ast$-algebra on the domain $\D=\rho(\B)\Hh$ such that 
$(\tilde{\rho}(\X), \rho(\B))$ is a compatible $O^\ast$-pair on the 
domain $\D$. That is, any continuous $\ast$-representation $\rho$ of $\B$ 
gives raise to a homomorphism of the (abstract) compatible pair $(\X,\B)$ 
to the compatible $O^\ast$-pair $(\tilde{\rho}(\X), \rho(\B))$.

\mn
{\bf Remark.} Let $\overline{\bbbb}(\D)$ be the completion of 
the normed $\ast$-algebra $(\bbbb(\D),\parallel{\cdot}\parallel)$. 
Obviously, the closure of the finite rank operators in 
$\overline{\bbbb}(\D)$ is 
the $\ast$-algebra $\C(\Hh)$ of compact operators on $\Hh$. Thus 
$\C(\Hh)$ is contained in the $C^\ast$-algebra $\overline{\bbbb}(\D)$. 
We call the quotient $C^\ast$-algebra
$$
C^\ast(\D):=\overline{B} (\D)/\C(\Hh)
$$
the $C^\ast$-{\it algebra associated with the domain} $\D$. It carries 
important information about the infinite dimensional closed subspaces of 
$\Hh$ contained in $\D$. As a sample, we mention the following result 
which is stated here without proof:
\\Suppose that $\D$ is a commutatively dominated Frechet 
domain (see [S1], p. 108, for the definition). Then the domain $\D$ 
contains an {\it infinite dimensional} closed linear subspace of $\Hh$ if and 
only if $C^\ast(\D)$ is non-trivial, that is, $C^\ast(\D)\ne\{0\}$.

\mn
{\bf 2. Well-behaved Representations of the Polynomial Algebra 
${\bbbc[x_1,{\dots}, x_{n}]}$}

\sn
In this section $\X$ denotes the $\ast$-algebra $\bbbc[x_1,{\dots}, x_n]$ of 
all polynomials with complex coefficients in $n$ commuting hermitean 
indeterminants $x_1,{\dots}, x_n$. 

\mn
{\bf Example 4:} Let $M$ be a closed subset of $\bbbr^n$ and let $\B_1$ be 
the $\ast$-algebra $C_0(\bbbr^n)$ 
of all compactly supported continuous functions on $\bbbr^n$  
with pointwise multiplication $(fg)(t)=f(t)g(t)$ and involution 
$f^+ (t)=\overline{f(t)}$. Let $\parallel f \parallel$ be the supremum 
of the function $|f(t)|$ over $M$. It is obvious that the multiplication
\begin{equation}\label{polact}
(p\ang f)(t_1,{\dots}, t_n)=p(t_1,{\dots}, t_n)f(t_1,{\dots}, t_n),~ 
p\in\X, f\in\B_1
\end{equation}
defines a left action of $\X$ on $\B_1$ such that $(\X,\B_1)$ is a compatible 
pair.

Now let $\rho$ be a non-degenerate continuous $\ast$-representation of 
$\B_1$ on a Hilbert space $\Hh$. It is well-known that there exists a 
spectral measure $E(\lambda),\lambda\in\bbbr^n$, on $\Hh$ supported in $M$ 
such that $\rho(f)=\int f(\lambda) dE(\lambda), f\in\B_1$. Then
$$
A_1:=\int\lambda_1 dE(\lambda),{\dots}, A_n:=\int\lambda_n dE(\lambda)
$$
are self-adjoint operators with commuting spectral projections and
\begin{align*}
\rho^\prime (x_j)\rho (f)=\rho (x_j\ang f)= \int\lambda_j f(\lambda) 
dE(\lambda)=\int\lambda_j dE(\lambda)\int f(\lambda) dE(\lambda)
= A_j \rho(f)
\end{align*}
for $j=1,{\dots}, n$. Conversely, any spectral measure on $\Hh$ with support 
contained in $M$ gives a $\ast$-representation $\rho^\prime$ as above.

We now specialize to the case where $M=\bbbr^n$. It is obvious that the 
operator 
$\rho^\prime(p)$ 
is essentially self-adjoint on $\rho(\B_1)\Hh$ for any $p=p^+\in\X$. 
Therefore, 
the $\ast$-representation $\rho^\prime$ (which is by definition the 
closure of its restriction to $\rho (\B_1)\Hh$ is integrable in the 
sense of [S], Chapter 9 (see Theorem 9.12). 
Conversely, 
any integrable $\ast$-representation of $\X=\bbbc[x_1,{\dots},x_n]$ is of 
this form. Thus, the $\ast$-representations $\rho^\prime$ associated with the 
compatible pair $(\X,\B_1)$ for $M=\bbbr^n$ are precisely the 
integrable $\ast$-representations of the polynomial algebra 
$\bbbc[x_1,{\dots},x_n]$.

\mn
{\bf Example 5.} Suppose that $K$ is a fixed compact subset $\bbbr^n$.
Let $\B_2=\X=\bbbc[x_1,{\dots},x_n]$ equipped with the supremum norm over 
the compact set $K$. With the left action (\ref{polact}) of $\X$ on $\B_2$, 
 $(\X,\B)$ is a compatible pair. 

Let $\rho$ be a continuous $\ast$-representation of $\B_2$ on a Hilbert 
space $\Hh$. Since $\rho$ is $\parallel{\cdot}\parallel$-continuous, 
$\rho$ extends to a $\ast$-representation of the $C^\ast$-algebra
$C(K)$. Hence there exists a spectral measure $E$ an $\Hh$ 
supported on the set $K$ such that $\rho (p)=\int p(\lambda) dE(\lambda)$ 
for $p\in\B_2$. As in the preceding example, we obtain
$$
\overline{\rho\prime (p)} 
=\int_K p(\lambda) dE(\lambda), ~p\in\X.
$$
Thus, the $\ast$-representations $\rho^\prime$ of $\X$ are precisely 
those bounded $\ast$-representation of $\bbbc[x_1\dots,x_n]$ for which the 
joint spectrum of the self-adjoint operators 
$\overline{\rho^\prime(x_1)},{\dots},\overline{\rho^\prime(x_n)}$ is 
contained in the set $K$. Among others, this example shows that the class of 
$\ast$-representations $\rho^\prime$ essentially depends on the choice 
of the norm $\parallel \cdot \parallel$.

\bn
{\bf 3. Integrable Representations of Enveloping Algebras}

\mn
Throughout this example, $G$ is a finite dimensional real Lie group with 
left Haar measure $\mu_l$ and Lie algebra $\g$ and $\E(\g)$ is the complex 
universal eneveloping algebra of $\g$. 

The algebra $\E(\g)$ is a $\ast$-algebra with involution determined by 
$x^+=-x$ for $x\in\g$. Let $\X$ denote the $\ast$-algebra 
$\E(\g)$. The vector space $B=C^\infty_0(G)$ is a 
$\ast$-algebra with respect to the convolution multiplication
\begin{equation}\label{conv}
(a\cdot b)(g)={\int\nolimits_G} a(h)b(h^{-1}g)d_{\mu_l}(h),
\quad a,b\in C^\infty_0(G),
\end{equation}
and the involution
\begin{equation}\label{inv}
a^+(g)=m(g)^{-1} \overline{a(g^{-1})}, a\in C^\infty_0(G),
\end{equation}
where $m$ denotes the modular function of the Lie group $G$. Since $m$ is 
a $C^\infty$-function on $G$ (see [Wa]), $a^+$ is again in $C^\infty_0(G)$. 
We equipp the $\ast$-algebra $\B$ with the $\ast$-invariant 
submultiplicative norm
$$
\parallel a\parallel=\int_G|a(g)|d\mu_l(g).
$$
All these facts are well-known and can be found, for instance, in
[Na], \S 28. The completion of $(\B,{\parallel\cdot\parallel})$ is 
nothing but the Banach $\ast$-algebra $L^1(G)$. 

Now we define the left action of $\X$ on $\B$. Let $x\rightarrow e^x$ 
denote the exponential map of $\g$ into $G$. Each element $x$ of 
$\E(\g)$ acts as a right-invariant differential operator $\tilde{x}$ of 
$G$. For $x\in\g$, the operator $\tilde{x}$ is given by
\begin{equation}\label{right}
(\tilde{x} a) (g):={\tfrac{d}{dt}}\big|_{t=0}~ a(e^{-tx} g),\quad a\in 
C^\infty_0(G).
\end{equation}
Using the formulas (\ref{conv}) and (\ref{right}) and the left 
invariance of the measure $\mu$ one easily verifies that 
\begin{equation}\label{left}
\tilde{x} (a{\cdot} b)=(\tilde{x} a){\cdot} b,~ a,b\in C^\infty_0(G),
\end{equation}
for $x\in \g$. Since the map $x\rightarrow \tilde{x}$ of $\E(\g)$ into the 
differential operators on $G$ is an algebra homomorphism, (\ref{left}) is valid for all $x\in\E(\g)$. {}From the preceding we conclude at once that
\begin{equation}\label{act}
x\ang a :=\tilde{x} a,\quad a\in C^\infty_0(G), x\in\E(\g),
\end{equation}
defines a left action of the $\ast$-algebra $\X=\E(\g)$ on the $\ast$-algebra $\B=C^\infty_0(G)$.

\mn
{\bf Lemma 6.} {\it $(\X, \B)$ is a compatible pair.}

\mn
{\bf Proof.} Clearly, it suffices to verify the compatibility 
condition (\ref{cond}) for elements $x$ of the Lie algebra $\g$. {}From 
the analysis on locally compact groups it is well-known
 (see [Na],  p. 376) that there is a {\it right} Haar measure $\mu_r$ on $G$ such that
\begin{equation}\label{rh}
d\mu_l(g)=m(g) d\mu_r(g),~ g\in G.
\end{equation}
Using formulas (\ref{conv}),(\ref{inv}),(\ref{right}),(\ref{act}) and (\ref{rh}) and 
\begin{align*}
(b^+{\cdot}(x^+\ang a))(g)&=\int b^+(h)(\tilde{x} a)(h^{-1} g) d\mu_l(h)\\
&=-{\tfrac{d}{dt}}\Big|_{t=0} \int m(h)^{-1} b(h^{-1}) a(e^{-tx} h^{-1}g) d\mu_l(h)\\
&=-{\tfrac{d}{dt}}\Big|_{t=0}\int b(h^{-1})a(e^{-tx} h^{-1}g) d\mu_r(h)\\
&=-{\tfrac{d}{dt}}\Big|_{t=0}\int b(e^{tx}k^{-1})a(k^{-1}g) d\mu_r(k e^{.tx})\\
&=-{\tfrac{d}{dt}}\Big|_{t=0}\int b(e^{tx}k^{-1})a(k^{-1}g) d\mu_r(k)\\
&=\int (\tilde{x}b)(k^{-1}) a(k^{-1}g) d\mu_r(k)\\
&=\int m(k^{-1}(\tilde{x} b)(k^{-1}a(k^{-1}g) d\mu_l(k)\\
&=\int (x\ang b)^+(k) a(k^{-1} g) d\mu_l(k)\\
&=((x\ang b)^+{\cdot}a)(g)
\end{align*}
for $x\in \g, a,b\in C^\infty_0(G)$ and $g\in G$. This proves that 
condition (\ref{cond}) is satisfied.\hfill $\Box$\\

Next let us look at the corresponding $\ast$-representations. In order to
 do so, we need to recall some facts on representation theory of Lie 
groups and Lie algebras which can be found in [S1], Chapter 10 or in [Wa]. Let $U$ be a 
unitary representation of the Lie group $G$ on a Hilbert space $\Hh$, 
that is, $g\rightarrow U(g)$ is a homomorphism of $G$ into the group of 
unitaries of $\Hh$ such that the map $g\rightarrow U(g)\varphi$ of $G$ 
into $\Hh$ is continuous for each veactor $\varphi\in\Hh$. Then there 
exists a unique $\ast$-representation $dU$ of the $\ast$-algebra $\E(\g)$ 
on the domain $\D^\infty(U)$ of $C^\infty$-vectors for $U$. For $x\in\g$, 
the operator $dU(x)$ acts as 
\begin{equation}\label{dax}
dU(x)\varphi={\tfrac{d}{dt}}\big|_{t=0} U(e^{tx})\varphi, \varphi\in\D^\infty(U).
\end{equation}
The linear span $\D_G(U)$ of vectors
$$
U_a\varphi:=\int a(g)U(g)\varphi~ d\mu_l(g),~ a\in C^\infty_0(G),
\varphi\in\Hh,
$$
is contained in the space $\D^\infty (U)$ of $C^\infty$-vectors. 
The vector space $\D_G(U)$ is called the G\.arding space of the 
unitary representation $U$. It was proved in [DM] that the G\.arding 
space $\D_G(U)$ is equal to $\D^\infty (U)$, but we shall not need this 
deep result here. For our purposes it is sufficient to know (see [S1], 
Corollary 10.1.16) that the G\.arding space is a core for all operators 
$dU(x),x{\in}\E(\g)$. This implies that the $\ast$-representation $dU$ on 
the domain $\D^\infty(U)$ is the closure of its restriction to $\D_G(U)$.

Suppose now that $\rho$ is a non-degenerate 
$\parallel\cdot\parallel$-continuous $\ast$-representation of 
the $\ast$-algebra $\B$ on the Hilbert space $\Hh$. Note that 
 all operators $\rho(x), x\in\B$, are bounded and defined on the 
 whole Hilbert space $\Hh$. Since $\B$ is $\parallel\cdot\parallel$-dense 
 in $L^1(G)$ and $\rho$ is $\parallel\cdot\parallel$-continuous, 
 $\rho$ extends by continuity to a non-degenerate $\ast$-representation, 
 denoted again by $\rho$, of the Banach $\ast$-algebra $L^1(G)$. It is 
 well-known ([Na], \sS 29, Theorem 1) that there exists a unique unitary 
 representaiton $U$ of the Lie group $G$ on the Hilbert space $\Hh$ such 
 that
\begin{equation}\label{uapi}
U_a=\rho(a), a\in L^1(G).
\end{equation}
Using formulas (\ref{dax}), (\ref{act}) and (\ref{uapi}), we obtain that
\begin{equation}\label{pi}
\rho^\prime(x) U_a=\rho^\prime(x)\rho (a)=\rho(x\ang a)=\rho(\tilde{x} a)= U_{\tilde{x} a}=dU(x)U_a
\end{equation}
for $x\in\E(\g)$ and $a\in C^\infty_0(G)$, where used the 
relation $dU(x)U_a=U_{\tilde{x}a}$ ([S1], Lemma 10.1.12). Therefore, 
the $\ast$-representation $\rho^\prime$ and $dU$ of $\X=\E(\g)$ coincide 
on the G\.arding space. Since $\rho^\prime$ and $dU$ are both the 
closures of their restriction to $\rho (\B)\Hh=\D_G(U)$, we conclude that 
$\rho^\prime=dU$.

Conversely, for any unitary representation $U$ there exists a unique 
non-degenerate $\ast$-representation $\rho$ of $L^1(G)$ and so of 
$\B=C^\infty_0(G)$ such that (\ref{uapi}) holds. By the above 
reasoning, we then have $\rho^\prime=dU$.

Summarizing, we have shown that the $\ast$-representations of the 
$\ast$-algebra $\X$ derived from the pair $(\X,\B)$ are precisely the 
$G$-integrable representations of the $\ast$-algebra $\E(\g)$ 
(in the sense of [S1], Chapter 10). That is, the $\ast$-representations 
$\rho^\prime$ are the $\ast$-representations $dU$ for 
unitary representations  
$U$ of the Lie group $G$.

\bn
{\bf 4. Examples related to the Weyl Calculus}

\mn
The Weyl calculus of pseudodifferential operators on $\bbbr^n$ can be used 
to construct further examples of compatible pairs. 
We restrict ourselves to the case $n=1$ and refer to the books 
[Fo] and [St] (see also [GV]) for the notation and the facts on the 
Weyl calculus needed in what follows.

Let $\cP$ and $\Q$ be the self-adjoint operators and let $W(s,t)$ be the 
unitary operator on the Hilbert space $L^2(\bbbr)$ defined by
\begin{align*}
(\cP f)(x)={\tfrac{1}{2\pi i}} ~f^\prime(x),~~ (\Q f)(x)=xf(x),~~
W(s,t)=e^{2\pi i(s\Q+t\cP)}, s,t\in\bbbr.
\end{align*}
To any measurable function $a$ on $\bbbr^2$ such that its 
Fourier transform
\begin{equation}\label{ahat}
\hat{a}(x,y)=\iint e^{-2\pi i(xs+yt)} a(s,t)~ dsdt
\end{equation}
is in $L^1(\bbbr^2)$, the Weyl calculus assigns an operator $Op(a)$ 
on the Hilbert space $L^2(\bbbr)$ by
\begin{equation}
Op(a)=\iint \hat{a}(s,t) W(s,t)~ dsdt.
\end{equation}
The integral is defined as a Bochner integral because 
$\hat{a}\in L^1(\bbbr^2)$. The adjoint operator $Op(a)^\ast$ and 
the operator product $Op(a) Op(b)$ (at least for "nice" 
symbols $a$ and $b$) are given by
\begin{equation}\label{oppro}
Op(a)^\ast=Op(a^+)~{\rm and}~ Op(a) Op(b)=Op(a\# b),
\end{equation}
where 
\begin{equation}\label{invol}
a^+(x,y):= \overline{a(x,y)},
\end{equation}
\begin{equation}\label{tp}
(a\# b)(x_1,x_2):=\iiiint a(u_1,u_2) b(v_1,v_2)e^{4\pi i[(x_1-u_1)(x_2-v_2)-(x_1-v_1)(x_2-u_2)]} du_1 du_2 dv_1 dv_2.
\end{equation}

Because of the formulas (\ref{oppro}) it is natural to expect that a 
vector space $\B$ of sufficiently nice functions becomes a $\ast$-algebra 
with involution (\ref{invol}) and product (\ref{tp}) provided 
that a $\in\B$ and $a\# b\in\B$ when $a,b\in\B$. An example of such a 
$\ast$-algebra is the Schwartz space $\sS(\bbbr^2)$. Another example is 
obtained as follows: Let $\A(\bbbr^2)$ denote the vector space of all 
holomorphic functions $f$ on $\bbbc^2$ such that for all 
$s_j, c_j,d_j\in\bbbr,c_j<d_j,j=1,2$, we have 
$$
\sup {\int\limits^\infty_{-\infty}}{\int\limits^\infty_{-\infty}} 
\left|f(x_1{+}iy_1, x_2{+}iy_2)\right|^2e^{s_1x_1+s_2x_2} dx_1dx_2<\infty,
$$
where the supremum is taken over the set
$\{(y_1,y_2)\in\bbbr^2: c_j<y_j<d_j,j=1,2\}$. Then $\B_1:=\sS(\bbbr^2)$ 
and $\B_2:=\A(\bbbr^2)$ are both $\ast$-algebras with product and 
involution defined by (\ref{tp}) and (\ref{invol}), respectively. 
This was shown in [GV], Proposition 1, for $\B_1:=\sS(\bbbr^2)$ and 
in [S4], Corollary 10, for $B_2=\A(\bbbr^2)$.

For $a \in \B_j$, $j=1,2$, the operator $Op(a)$ is bounded. 
{}From the formulas (\ref{oppro}) it follows that the norm 
$\parallel \cdot \parallel$ on $\B_j$ defined by 
$$
\parallel a \parallel:=\parallel Op(a) \parallel
$$
is a submultiplicative and $\ast$-invariant. Therefore, 
$\B_1$ and $\B_2$ are normed $\ast$-algebras.

Let $\X_1$ be the unital algebra with two generators $\pb$ and $\bx$ and defining relation
$$
\pb\bx -\bx\bf{p}=-i1.
$$
Clearly, $\X_1$ is a $\ast$-algebra with involution determined by 
$\bf{p}^+=\bf{p}$ and $\bx^+=\bx$. One easily checks that there is a 
left action of the algebra $\X_1$ on 
$\B_1=\sS(\bbbr^2)$ 
such that
\begin{equation}\label{actone}
\pb\ang a:=\left( {1\over{2i}}~{\partial\over{\partial x_1}}+2\pi x_2\right) a,~\bx\ang a:=\left(x_1-{1\over{4\pi i}}{\partial\over{\partial x_2}}\right) a
\end{equation}
In terms of the operators $\cP_j={1\over{2\pi_i}} 
{\partial\over{\partial x_j}}$ and $Q_j=x_j$ on the 
Hilbert space $L^2(\bbbr^2)$ the latter can be rewritten as 
\begin{equation}\label{actonez}
\bf{p}=\pi\cP_1+2\pi Q_2,~ x=Q_1-{1\over 2}\cP_2.
\end{equation}

Suppose that $q$ is a complex number of modulus one. Let 
$\X_2$ be the coordinate $\ast$-algebra $\cO(\bbbr^2_q)$ of the real 
quantum plane. It is defined as follows (see, for instance, [FRT] or [KS]). 
As an algebra, 
$\cO(\bbbr^2_q)$ has two generators $\bx$ and $\by$ with defining relation
$$
\mathbf{x}\mathbf{y}= q\mathbf{y}\mathbf{x}
$$
The involution is defined by the requirements $\bx^+=\bx$ and 
$\by^+=\by$. 

We write $q=e^{2\pi i\gamma}$ with a {\it fixed} real number $\gamma$ and 
we take two real numbers $\alpha$ and $\beta$ such that 
$\alpha\beta=\gamma$. Then we define a left action of the algebra $\X_2$ on $\B_2=\A(\bbbr^2)$ by 
\begin{equation}\label{acttwo}
(\bx\ang a)(x_1,x_2)=e^{2\pi\alpha x_1} a(x_1, x_2+i\alpha /2), 
(\by\ang a)(x_1,x_2)=e^{2\pi\beta x_2} a( x_1-i\beta /2, x_2).
\end{equation}
Since $x\ang (y\ang a)=qy\ang(x\ang a)$ as easily verified, the latter gives a well-defined left action of $\X_2$ on $\B_2$.
The operator $e^{2\pi c\cP},c\in\bbbr$, acts on functions of its 
domain as 
$$
\left( e^{2\pi c \cP} f\right)(x)=f(x-ci);
$$ 
see Lemma 1.1 in [S3] for a precise statement. 
Using this fact, formula (\ref{acttwo}) means that 
\begin{equation}\label{acttwoz}
\bx=e^{2\pi\alpha\Q_1}\otimes e^{-\pi{\alpha\over 2} \cP_2},
~\by=e^{\pi{\beta\over 2}\cP_1} \otimes e^{2\pi\beta\Q_2}.
\end{equation}

\bn
{\bf Lemma 7.} {\it $(\X_1, \B_1)$ and $(\X_2, \B_2)$ are 
compatible pairs.}

\mn
{\bf Proof.} It suffices to check (\ref{cond}) for the 
generators  $\bf{p},\bx$ and $\bx,\by$, respectively. By (\ref{actonez}), 
for the generators $\bf{p}$ and $\bx$ condition (\ref{cond}) means that 
$$
(\cP_1a+2\Q_2 a)^+\# b=a^+\#(\cP_1b+2\Q_2 b),(2\Q_1a-\cP_2a)^+\# b=a^+\# (2\Q_1 b-\cP_2 b).
$$
Both relations are easily derived from the definition (\ref{tp}) of the 
twisted product $\#$ using partial integration. Note that the 
corresponding boundary terms vanish because the functions $a$ and  $b$ are 
in $\sS(\bbbr^2)$. We omit the details.

For the generators $\bx$ and $\by$ of $\X_2$ formula (\ref{cond}) says that 
\begin{align*}
&\left( e^{2\pi\alpha\Q_1} e^{-\pi{\alpha\over 2}\cP_2} 
a\right)^+\# b=a^+\#\left( e^{2\pi\alpha\Q_1} e^{-\pi{\alpha\over 2}\cP} b\right),\\
&\left( e^{\pi{\beta\over 2}\cP_1} e^{2\pi\beta\Q_2} a\right)^+\# 
b=a^+\#\left(e^{\pi{\beta\over 2}\cP_1} e^{2\pi\beta\Q_2} b\right).
\end{align*}
Both identities follow at once from the formulas stated in [S4], Lemma 9, 
combined with the fact that $\left(e^{c\cP_j}a\right)^+=e^{-c\cP_j} a^+$ for 
$c\in\bbbr$ and  $a\in\A(\bbbr^2)$.\hfill$\Box$

\sn
Now we turn to $\ast$-representations. If $\K$ is a Hilbert space, 
then it obvious from (\ref{oppro}) that the formula
\begin{equation}\label{rep}
\rho_0(a)= Op(a)\otimes I,~ a\in\B_j~,
\end{equation}
defines a continuous $\ast$-representation of the normed $\ast$-algebra 
$\B_j$ on the Hilbert space $L^2(\bbbr)\otimes \K$.

\mn{\bf Lemma 8.} {\it Any continuous $\ast$-representation of the normed 
$\ast$-algebra $\B_j,j=1,2$, is unitarily equivalent to a 
$\ast$-representation $\rho_0$.}

\sn
This assertion is probably well-known, but we could not find it 
explicitely in the literature. Thus we include a sketch of proof. 

\mn
{\bf Sketch of proof.} Suppose that $\rho$ is a continuous 
$\ast$-representation of the normed $\ast$-algebra $\B_j$ on a 
Hilbert space $\G$. Since $\rho$ is a direct sum of cyclic 
$\ast$-representation, we can assume without loss of generality 
that $\rho$ is cyclic. Let $\varphi\in\G$ be a cyclic vector for $\rho$. 
For $a\in\B_j$ and $s,t\in\bbbr$, we set 
\begin{equation}\label{ast}
a_{s,0}(x_1,x_2):=e^{2\pi isx_1} a(x_1,x_2{-}s/2), 
a_{0,t}(x_1,x_2):=e^{2\pi tx_2} a(x_1{+}t/2,x_2).
\end{equation}
Since $a_{s,0}^+\# a_{s,0}=a^+\# a$ as easily computed, we have 

\begin{align*}
\parallel\rho(a_{s,0})\varphi\parallel^2&=
\langle\rho(a_{s,0})\varphi,\rho(a_{s,0})\varphi\rangle=
\langle\rho(a_{s,0}^+\# a_{s,0})\varphi,\varphi\rangle\\
&=\langle\rho(a^+\# a) \varphi,\varphi\rangle=
\langle\rho(a)\varphi,\rho(a)\varphi\rangle=
\parallel\varphi (a)\rho\parallel^2.
\end{align*}
Hence there is an isometric map $U(s)$ of $\rho(\B_j)\varphi$ onto 
itself such that $U(s)\rho(a)\varphi=\rho(a_{s,0})\varphi,a\in\B_j$. 
Obviously, $U(s_1{+}s_2)=U(s_1)U(s_2)$ for $s_1,s_2\in\bbbr$ and $U(0)=I$. 
Moreover, $\parallel a_{s,0}-a\parallel_0\rightarrow 0$ as 
$s\rightarrow 0$ and hence $U(s)\psi\rightarrow\psi$ in $\G$ 
for $\psi\in\rho(\B_j)\varphi$ as $s\rightarrow 0$ by the continuity of 
$\rho$. Since $\rho(\B_j)\varphi$ is dense in $\G$, $U(s)$ extends to a 
unitary operator on $\G$ and $s\rightarrow U(s)$ is a strongly 
continuous one-parameter unitary group on $\G$. Similarly, there is 
another strongly continuous one-parameter unitary group 
$t\rightarrow V(t)$ on $\G$ such that $V(t)\rho(a)\varphi=
\rho(a_{0,t})\varphi$. {}From their definitions we immediately 
derive that the unitary groups $U$ and $V$ satisfy 
the Weyl relation
$$
V(t)U(s)=e^{2\pi ist} U(s)V(t),s,t\in\bbbr.
$$
Therefore, by the Stone- von Neumann uniqueness theorem for the canonical 
commutation relation (see e.g. [Pu]), there exists a Hilbert space $\K$ 
and a unitary map $T$ of $\G$ onto $L^2(\bbbr)\otimes \K$ such that
\begin{equation}\label{uniqu}
T^{-1} U(s)T=W(s,0)\otimes I,~ T^{-1}V(t)T=W(0,t)\otimes 
I~{\rm for}~s,t\in\bbbr.
\end{equation}
Let us abbreviate $\tilde{W}(s,t):=e^{-\pi i st}U(s)V(t)$ and 
$\tilde{\rho}(a):=\iint \hat{a}(s,t)\tilde{W}(s,t) ds dt$, where 
$a\in\B_j$. Since
\begin{equation}\label{wrel}
W(s,t)=e^{\pi ist} W(s,0)W(0,t),s,t\in\bbbr,
\end{equation}
it follows from (\ref{uniqu}) that $T^{-1}\tilde{W}(s,t)T=W(s,t)\otimes I$ 
and hence
\begin{equation}\label{tpi}
T^{-1}\tilde{\rho} (a)T=Op(a)\otimes I=\rho_0(a),~ a\in\B_j.
\end{equation}
On the other hand from the definition of $Op(a)$ one derives that 
$$
W(s,0)Op(a)=Op(a_{s,0})~{\rm and}~ W(0,t)Op(a)=Op(a_{0,t}).
$$
By the defnition of $\tilde{W}(s,0)=U(s)$ and $\tilde{W}(0,t)=V(t)$, the 
latter implies that
$$
\tilde{W}(s,t)\rho(b)\varphi=\rho(Op^{-1}(W(s,t)Op(b)))\varphi
$$
and so
\begin{align*}
&\tilde{\rho}(a)\rho(b)\varphi=\rho (Op^{-1}\Big(\Big(\iint\hat{a}(s,t)
W(s,t)ds dt)Op(b)\Big)\Big)\varphi\\
&=\rho(Op^{-1}(Op(a)Op(b)))\varphi=\rho(Op^{-1}(Op(a\# b)))\varphi=
\rho(a)\rho(b)\varphi
\end{align*}
for $a,b\in\B_j$. Since $\rho(\B_j)\varphi$ is dense in $\G$, we obtain 
$\tilde{\rho}(a)=\rho(a)$. Thus, by (\ref{tpi}) we have  
$T^{-1}\rho(a)T=\rho_0(a)$ which completes the proof of Lemma 8.\hfill $\Box$

\sn
Finally, let us describe the $\ast$-representation $\rho_0^\prime$ of the 
$\ast$-algebra $\X_j$ derived from the $\ast$-representation $\rho_0$ of 
$\B_j$. First let $j=1$. Since $W(s,0)=e^{2\pi is\Q}$ and 
$W(0,t)=e^{2\pi it\cP}$, it follows from (\ref{ast}) and 
(\ref{wrel}) by differentation at $s=0$ and $t=0$, respectively, that
$$
\Q Op(a)=Op\Big(x_1a-{1\over{2\pi i}}{{\partial a}
\over{\partial x_2}}\Big),~ \cP Op(a)=Op\left({1\over{2i}}
{{\partial a}\over{\partial x_1}}+2\pi x_2a\right).
$$
Combining the latter with (\ref{actone}) we conclude that 
$$
\rho^\prime_0 (\bx)\rho_0(a)=\rho_0(\bx\ang a)=Op(\bx\ang a)\otimes 
I=\Q Op(a)\otimes I=(\Q\otimes I)\rho_0(a)
$$
and similarly $\rho^\prime_0({\bf{p}})\rho_0(a)=(\cP\otimes I)\rho_0(a)$ 
for $a\in\B_1$, so that 
$$
\overline{\rho^\prime_0(\bx)}=\Q\otimes I~{\rm and}~ 
\overline{\rho^\prime_0(\bf{p})}=\cP\otimes I.
$$
That is, up to unitary equivalence the $\ast$-representations 
$\rho^\prime_0$ of $\X_1$ are precisely the orthogonal direct sums of the 
Schr"odinger representation of the $\ast$-algebra $\X_1$ with 
domain $\sS(\bbbr)$ on the Hilbert space $L^2(\bbbr)$.

Now let $j=2$. {}From the formulas in [S4], Lemma 9, we then have 
\begin{equation}\label{opa}
e^{2\pi\alpha \Q} Op(a)=Op(\bx\ang a),~ e^{2\pi\beta\cP}Op(a) = 
Op(\by\ang a),
\end{equation}
where $x\ang a$ and $y\ang  a$ are defined by (\ref{acttwo}). 
{}From (\ref{opa}) and (\ref{def}) we obtain
$$
\overline{\rho^\prime_0(\bx)}=e^{2\pi\alpha\Q}\otimes I,
\quad \overline{\rho^\prime_0(\by)}=e^{2\pi\beta\cP}\otimes I.
$$
This $\ast$-representation $\rho_0$ appears (in a slighty different notation) 
in the work by M. Rieffel [R] on the quantum plane.
 
\sn
We illustrate the preceding considerations by three other closely related 
examples. Let $\B_3$ be the normed $\ast$-algebra $\B_2\oplus\B_2\oplus\B_2
\oplus\B_2$. There is a left action of $\X_2$ on $\B_3$ such that
\begin{align*}
&\bx\ang (a_1,a_2,a_3,a_4)=(\bx\ang a_1,\bx\ang a_2,
-\bx\ang a_3,-\bx\ang a_4),\\
&\by\ang (a_1,a_2,a_3,a_4)=
(\by\ang a_1,-\by\ang a_2, \by\ang a_3, -\by\ang a_4)
\end{align*}
for $a_1,a_2,a_3,a_4\in\B_2$, where $\bx\ang a$ and $\by\ang a$ 
are defined by (\ref{acttwoz}). Obviously, $(\X_2,\B_3)$ is a 
compatible pair.

For arbitrary $\epsilon_1,\epsilon_2\in\{+,-\}$, there is a 
$\ast$-representation $\rho^\prime_{\epsilon_1,\epsilon_2}$ of 
$\X_2$ derived from a continuous $\ast$-representation 
$\rho_{\epsilon_1,\epsilon_2}$ of $\B_3$ such that
$$
\overline{\rho^\prime_{\varepsilon_1,\varepsilon_2}(\bx)}=
\varepsilon_1 e^{2\pi\alpha\Q}\otimes I,
~\overline{\rho^\prime_{\varepsilon_1,\varepsilon_2}(\by)}=
\varepsilon_2 e^{2\pi\beta\cP}\otimes I.
$$
It is easily seen that the $\ast$-representation $\rho^\prime$ of $\X_2$ associated with the 
pair $(\X_2,\B_3)$ are precisely the orthogonal direct sums of some of 
$\ast$-representations $\rho^\prime_{++},\rho^\prime_{+,-},\rho^\prime_{-+},
\rho^\prime_{--}$.
\mn

Finally, let $\B_4$ the $\ast$-algebra 
$\B_2\otimes M_2(\bbbc)\cong M_2(\B_2)$ equipped with the $C^\ast$-matrix 
norm derived from the $C^\ast$-norm of $\B_2$. We suppose now that 
$\alpha$ and $\beta$ are real numbers such that
\begin{equation}\label{alphabeta}
\alpha\beta=\gamma+{1/2}.
\end{equation}
Then there exists a left action of the algebra $\X_2$ on $\B_4$ given by

\begin{equation*}
\bx\ang \left( \begin{matrix} a_1 &a_2\\ a_3 &a_4\end{matrix} \right) = 
       \left( \begin{matrix}~~\bx\ang a_1 &~~\bx\ang a_2\\ 
                      -\bx\ang a_2  & -\bx\ang a_4\end{matrix} \right),~
\by\ang \left( \begin{matrix} a_1 &a_2\\ a_3 &a_4\end{matrix} \right) = 
       \left( \begin{matrix} \by\ang a_3  &\by\ang a_4\\ 
                        \by\ang a_1 &\by\ang a_2\end{matrix}\right).
\end{equation*}
\\
where
$$
\bx\ang a:=-e^{2\pi\alpha\Q_1}\otimes e^{-\pi{\alpha\over 2} \cP_2} a,
~ \by\ang e^{\pi{\beta\over 2}\cP_1}\otimes e^{2\pi\beta\Q_2}a.
$$
Note that the latter formula coincides with (\ref{acttwoz}), but now we have 
assumed $\alpha\beta=\gamma+{1\over 2}$ rather than $\alpha\beta=\gamma$. 
Using the fact that $(\X_2,\B_2)$ is a compatible pair, a straightforward 
computation shows that $(\X_2,\B_4)$ is also a compatible pair.

For a Hilbert space $\G$, there is a continuous $\ast$-representation 
$\rho_0$ of $\B_4$ on the Hilbert space $L^2(\bbbr)\otimes\bbbc^2\otimes\G$ 
such that 
$$
\rho_0(a\otimes m)=Op(a)\otimes m\otimes I,~a\in\B_2,m\in M_2(\bbbc).
$$
{}From Lemma 8 it follows immediately that any continuous 
$\ast$-representation of $\B_4$ is unitarly equivalent to such a 
$\ast$-representation $\rho_0$. Using the facts that 
$$
Op(\bx\ang a)=e^{2\pi\alpha\Q} Op(a),~Op(\by\ang a)=e^{2\pi\beta\cP}Op(a),
$$
we derive that the closures of the representation operator
 $\rho^\prime_0(\bx)$ and $\rho^\prime_0(\by)$ are the self-adjoint 
 operator matrices
$$
\overline{\rho^\prime_0 (\bx)}=
\left(\begin{matrix}e^{2\pi\alpha\Q}\otimes I &0\cr 0 
&-e^{2\pi\alpha\Q}\otimes I \end
{matrix}\right),~ \overline{\rho^\prime_0 (\by)}=
\left(\begin{matrix} 0&e^{2\pi\beta\cP}\otimes I 
\cr e^{2\pi\beta\cP}\otimes I &0\end{matrix}\right)~.
$$
\mn

Finally, we consider a $\ast$-algebra and their representations in the 
above context which was used by S.L. Woronowicz in his approach to the 
quantum 
$ax+b$-group [W2]. Let $\X_3$ denote the $\ast$-algebra generated by three 
hermitean elements $x,y,\chi$ and defining relations
$$
xy=qyx,~x\chi=\chi x,~ y\chi=-\chi y,~\chi^2=1
$$
Then there is a left action of $\X_3$ on the $\ast$-algebra 
$\B_4=\B_2\otimes M_2(\bbbc)$ defined by

\begin{align*}
&x \ang \left( \begin{matrix} a_1 &a_2\\ a_3 &a_4\end{matrix}\right)
       =\left( \begin{matrix} x\ang a_1 &x\ang a_2\\  
                              x\ang a_3 &x\ang a_4\end{matrix}\right), ~
 y\ang \left(  \begin{matrix} a_1 &a_2\\  a_3 &a_4\end{matrix}\right) 
=      \left(  \begin{matrix} y\ang a_1 &y\ang a_2\\ 
                              -y\ang a_3 &-y\ang a_4\end{matrix}\right)~,\\
&\chi\ang\left(\begin{matrix} a_1 &a_2\\  a_3&a_4\end{matrix}\right)= 
       \left(  \begin{matrix} a_3 &a_4\\ a_1 &a_2\end{matrix}\right)~,
\end{align*}
where $x\ang a$ and $y\ang a$ are given by (\ref{acttwoz}). It is not 
difficult to check that $(\X_3,\B_4)$ is a compatible pair and that 
for the corresponding $\ast$-representation $\rho^\prime_0$ of $\X_3$ we 
have
$$
\overline{\rho^\prime_0(x)}=\left(\begin{matrix} e^{2\pi\alpha\Q}\otimes I 
&0\\ 0 &e^{2\pi\alpha\Q}\otimes I\end{matrix}\right),
~\overline{\rho^\prime_0(y)}=\left(\begin{matrix} 
e^{2\pi\beta\cP}\otimes I &0\\ 0 &-e^{2\pi\beta\cP} 
\otimes I\end{matrix}\right),
$$
where $\alpha\beta=\gamma$.

\bn
{\bf 5. Another  Example: Quantum  $SU_q(1,1)$ Group}

\mn

Suppose that $q$ is a real number such that $q\ne 0, 1,-1$. Let $\X$ denote 
the coordinate $\ast$-algebra $\cO(SU_q(1,1))$ of the quantum group 
$SU_q(1,1)$ (see, for instance, [FRT], [KS] or [W1]). That is, $\X$ is 
the $\ast$-algebra with two generators $a$ and $c$ and defining relations
$$
ac=q ca,~ ac^+ = qc^+ a,~ cc^+=c^+c, ~a^+ a-c^+ c=1,~aa^+ -q^2c^+ c=1.
$$
Let $\B$ be the $\ast$-algebra generated by the 
$\ast$-algebra $C_0(\bbbc)$ of compactly supported continuous 
functions on the complex plane $\bbbc$ with usual algebraic structure 
and a single generator $v$ satisfying the relations 
$(vf)(z)=f(q^{-1}z), f\in C_0(\bbbc), z\in\bbbc$, and $vv^+=v^+v=1$. 
We consider $\B$ as $\ast$-subalgebra of bounded operators acting on the 
Hilbert space $L^2(\bbbc)$ equipped with the operator norm.

Then there is a left action of $\X$ on $\B$ defined by
\begin{align*}
 (a\ang f)(z)&={\sqrt{1+|z|^2}}~ f(q^{-1}z),~(c\ang f)(z)=zf(z),\\
 a\ang v     &={\sqrt{1+|z|^2}}~ v^2,~c\ang v=zv.
\end{align*}
One easily checks that $(\X,\B)$ is a compatible pair.

Let $\rho$ be a continuous $\ast$-representation of $\B$ on a Hilbert 
space $\Hh$. It is not difficult to show that there exists a 
spectral measure $E(z), z\in\bbbc$, and a unitary operator $u$ on 
$\Hh$ satisfying $uE(z)u^{-1}=E(q z), z\in\bbbc$, such that
$$
\rho(f)=\int f(z) dE(z), f\in C_0(\bbbc),~{\rm and}~ \rho(v)=u.
$$
The operators $\overline{\rho^\prime(a)}$ and 
$\overline{\rho^\prime(c)}$ of the associated $\ast$-representation 
$\rho^\prime$ of $\X$ act as
$$ \rho^\prime(a)={\sqrt{1+(z)^2}}~ u~~{\rm and}~~\rho^\prime(c)=z.
$$
These are precisely the $\ast$-representations of the 
$\ast$-algebra $\X=\cO(SU_q(1,1))$ considered in [W1], p. 79.

\mn
{\bf 6. Outlook}

\sn
In the preceding sections we have investigated a variety of examples 
of compatible pairs $(\X,\B)$ and corresponding well-behaved 
$\ast$-representations. Here we add some remarks on these concepts 
and on possible modifications.

\sn
1. First, let us emphasize that the above approach does not solve the 
problem of selecting the well-behaved $\ast$- representations of a given 
$\ast$-algebra. As mentioned in the introduction, this essentially depends 
on the specific structure of the $\ast$-algebra and on the aim of the 
considerations. The notion of a compatible pair is only a proposal of a 
general concept in order to treat appearently different examples 
in the same general context.

2. For a compatible pair $(\X,\B)$, let $\tau$ denote the locally 
convex topology on $\B$ defined by the family of seminous 
$p_x(b)=\parallel x\ang b\parallel,x\in\X,b\in\B$. In all above 
examples, $\B[\tau]$ is a metrizable $\ast$-algebra with jointly 
continuous multiplication. Therefore, the completion 
$\tilde{\B}[\tau]$ of $\B[\tau]$ is a Frechet $\ast$-algebra and the 
action of $\X$ on $\B$ extends by continuity to an action on 
$\tilde{B}$ which also satisfies the compatibility condition
 (\ref{cond}). Thus, another possible approach might be to replace 
the (uncomplete) normed $\ast$-algebra $(\B,\parallel{\cdot}\parallel)$ 
by the Frechet $\ast$-algebra $\tilde{\B}[\tau]$.

3.  The reason why we have defined 
compatible pairs by means of condition (\ref{cond}) is that it 
seems to be the weakest requirement on the left action $\ang$ 
which ensures that a non-degenerate $\ast$-representation of $\B$ yields a 
well-defined $\ast$-representation of $\X$ by means of formula (\ref{def}). 
However, for compatible $O^\ast$-pairs (Example 3) and for all compatible 
pairs occuring in this paper we have much more structure: In all theses cases 
the $\ast$-algebras $\X$ and $\B$ 
are $\ast$-subalgebras of a larger $\ast$-algebra $\A$ and 
the left action $x \ang b$ is just the product $x{\cdot}b$ of the elements 
$x \in \X$ and $b \in \B$ in the algebra $\A$. It is obvious that the 
$\ast$-algebra axioms imply that condition (\ref{cond}) is fulfilled.

4. Let $\X$ be a $\ast$-algebra with unit and let $\ang$ be a left action 
of $\X$ on another $\ast$-algebra $\B$. On the direct sum $\A:= \X \oplus \B$ of 
vector spaces $\X$ and $\B$ we define a product 
$$
(x+a)(y+b):= xy+ (y^\ast \ang  a^\ast)^\ast + x áang b + ab,~ 
x,y \in \X, a,b \in \B
$$
and an involution $(x+a)^\ast:=x^\ast +a^\ast$. With these structures, $\A$ 
is a $\ast$-algebra if and only if conditions (\ref{cond}) and
$$
(x \ang a)b = x \ang (ab)~{\rm and}~ 
(x \ang(y\ang a)^\ast)^\ast = y \ang (x \ang a^\ast)^\ast,~ 
x,y \in \X, a,b \in \B,
$$
are fulfilled. If this is true, then 
$\X$ and $\B$ are $\ast$-subalgebras of $\A$ and the left action $\ang$ of $\X$ 
on $\B$ is the product in the algebra $\A$. Such a $\ast$-algebra $\A$ 
has been essentially used as ''function algebra'' on 
the quantum quarter plane in [S4].

5. The elements $x \in \X$ for a compatible pair$(\X,\B)$ can be also considered 
as multipliers of the algebra $\B$. 


\sn


\end{document}